\documentclass[pdflatex,sn-mathphys-num]{sn-jnl}

\usepackage{amsmath,amssymb,amsthm}

\theoremstyle{thmstyleone}
\newtheorem{theorem}{Theorem}
\newtheorem{lemma}[theorem]{Lemma}
\newtheorem{corollary}[theorem]{Corollary}
\newtheorem{proposition}[theorem]{Proposition}

\DeclareMathOperator{\dens}{dens}

\hypersetup{
  pdftitle={The modulus of smoothness and lp-generation},
  pdfauthor={Zhiyao Fang, Wenwen Zhang, and Qi Liu},
  pdfsubject={Power-type smoothness and Banach space generation},
  pdfkeywords={Banach space generation, modulus of smoothness, superreflexivity}
}

\renewcommand{\abstracthead}[2]{\par\unboldmath\noindent\textsc{#2.}\enspace}

\begin{document}

\title[Modulus of smoothness and \texorpdfstring{\(\ell_p\)}{lp}-generation]%
{The modulus of smoothness and \texorpdfstring{\(\ell_p\)}{lp}-generation}

\author[1]{\fnm{Zhiyao} \sur{Fang}}
\email{fzy73@stu.aqnu.edu.cn}

\author[1]{\fnm{Wenwen} \sur{Zhang}}
\email{Y26220004@stu.aqnu.edu.cn}

\author[1]{\fnm{Ying} \sur{Xu}}
\email{Y26220046@stu.aqnu.edu.cn}

\author*[1]{\fnm{Qi} \sur{Liu}}
\email{liuq67@aqnu.edu.cn}

\affil[1]{\orgdiv{School of Mathematics and Statistics},
  \orgname{Anqing Normal University},
  \orgaddress{\city{Anqing}, \postcode{246133}, \state{Anhui}, \country{China}}}

\abstract{Fabian, Godefroy, H\'ajek and Zizler
[\emph{J. Funct. Anal.} \textbf{200} (2003), 301--323, Theorem~2]
proved that at density \(\aleph_1\), generation by a space whose
modulus of smoothness has power type \(p\), \(1<p\le 2\), is
equivalent to \(\ell_p(\Gamma)\)-generation. We extend this
equivalence to arbitrary density, with \(|\Gamma|=\dens X\),
answering \mbox{Problem~291} of Guirao, Montesinos and Zizler
[\emph{Open Problems in the Geometry and Analysis of Banach Spaces},
Springer, Cham (2016), p.~127]. The proof combines an upper
\(\ell_p\)-estimate on successive projection differences
with transfinite induction on density. We characterize generation
by superreflexive spaces and by their \(\ell_2\)-sums, and obtain
\(\ell_{\min\{r,2\}}\)-generation for closed subspaces of
\(L^r(\mu)\), \(1<r<\infty\), including Hilbert generation for
\(r\ge 2\).}

\keywords{Banach space generation, Modulus of smoothness,
Projectional resolution of the identity, Superreflexivity}

\pacs[MSC Classification]{46B26, 46B20, 46B03}

\maketitle

\section{Introduction}\label{sec:introduction}

A Banach space \(X\) is \emph{generated by} a Banach space \(Y\) if
there is a bounded linear operator \(T:Y\to X\) whose range is norm
dense in \(X\). Generation by a Hilbert space is the case of
\(\ell_2(\Gamma)\)-generation. We write \(\dens X\) for the least
cardinality of a norm-dense subset of \(X\), and \(S_Y\) for the unit
sphere of \(Y\). The modulus of smoothness of a nonzero Banach space
\(Y\) is
\[
  \rho_Y(t)=
  \sup_{u,v\in S_Y}\left(\frac{\|u+tv\|+\|u-tv\|}{2}-1\right),
  \qquad t\ge 0 .
\]
It has power type \(p\) if \(\rho_Y(t)=O(t^p)\) as \(t\downarrow 0\).
We use the equivalent estimate \(\rho_Y(t)\le C t^p\) for \(0<t\le 1\),
which follows by enlarging \(C\) and using \(\rho_Y(t)\le t\).

Fabian, Godefroy, H\'ajek and Zizler established that, for spaces of
density \(\aleph_1\) and \(1<p\le 2\), generation by a Banach space
with a modulus of smoothness of power type \(p\) is equivalent to
\(\ell_p(\Gamma)\)-generation \cite[Theorem~2]{FGHZ}. The cardinal
restriction for this implication was explicitly raised in
\cite[Theorem~9 and the discussion thereafter]{FGMZ} and recorded as
Problem~291 in \cite[p.~127]{GMZ}. Related results of Fabian, Gonz\'alez
and Zizler concern generation by asymptotically \(p\)-flat sets at
density \(\aleph_1\) \cite[Theorem~5]{FGZ-flat}. At the quadratic
endpoint, the survey \cite[Theorem~24]{FMZ-survey} states a
Hilbert-generation criterion in terms of smoothness of order two
along a bounded set with dense linear span.

\begin{theorem}\label{thm:generation}
Let \(1<p\le 2\) and let \(X\) be a nonzero Banach space. The following
statements are equivalent.
\begin{enumerate}
\renewcommand{\labelenumi}{(\roman{enumi})}
\item \(X\) is generated by a Banach space admitting an equivalent
      norm whose modulus of smoothness has power type \(p\).
\item \(X\) is generated by \(\ell_p(\Gamma)\) for an index set
      \(\Gamma\) with \(|\Gamma|=\dens X\).
\end{enumerate}
\end{theorem}

The proof applies the projectional resolution to the \emph{generating}
space. Its difference spaces inherit the modulus of smoothness and
have smaller density. An upper \(\ell_p\)-estimate on finite sums of
vectors from these spaces permits their generating maps to be
assembled into a single map. The constructions work over either
scalar field.

The Day--Nordlander bound gives
\(\rho_Y(t)\ge \sqrt{1+t^2}-1\) when the underlying real dimension
of \(Y\) is at least two. This follows from the convexity bound in
\cite{Nordlander} and the duality of the moduli
\cite[Chapter~11]{Pisier}. Thus \(p>2\) in Problem~291 reduces to
the case of a generator of dimension at most one, for which the
conclusion follows from finite dimensionality.

\section{Contractive projections and upper estimates}\label{sec:projections}

\begin{lemma}\label{lem:one-sided}
Let \(1<p\le 2\), and suppose that the modulus of smoothness of
\(Y\) satisfies \(\rho_Y(t)\le C t^p\) for \(0<t\le 1\).
If \(P:Y\to Y\) is a contractive projection, \(u\in P Y\), and
\(v\in\ker P\), then
\begin{equation}\label{eq:one-sided}
  \|u+v\|^p \le \|u\|^p + D\|v\|^p,
  \qquad D=2^p\max\{1,pC\}.
\end{equation}
\end{lemma}

\begin{proof}
Contractivity gives \(\|u\|\le\|u+v\|\) and
\(\|u\|\le\|u-v\|\). Suppose first that \(0<\|v\|\le\|u\|\).
The definition of the modulus gives
\[
  \|u+v\|-\|u\|
  \le 2\|u\|\rho_Y(\|v\|/\|u\|)
  \le 2C\|u\|^{1-p}\|v\|^p .
\]
Since \(\|u+v\|\le 2\|u\|\), the mean value theorem applied to
\(s\mapsto s^p\) bounds \(\|u+v\|^p-\|u\|^p\) by
\(2^p pC\|v\|^p\). If \(\|v\|>\|u\|\), then
\(\|u+v\|^p\le 2^p\|v\|^p\). The case \(v=0\) is immediate.
\end{proof}

We use the following form of the classical projectional resolution
theorem \cite{AmirLindenstrauss}. We include the construction for
reflexive spaces to make the density bound at limit ordinals explicit.

\begin{lemma}\label{lem:pri}
Let \(Y\) be a nonseparable reflexive Banach space, and let \(\kappa\)
be the initial ordinal of cardinality \(\dens Y\). There are contractive
projections \(P_\alpha\) for \(0\le\alpha\le\kappa\), with
\(P_0=0\) and \(P_\kappa=I_Y\), such that
\[
  P_\alpha P_\beta=P_\beta P_\alpha=P_\alpha
  \quad(\alpha\le\beta),\qquad
  \dens(P_\alpha Y)\le \max\{\aleph_0,|\alpha|\}
  \quad(\alpha<\kappa).
\]
At each nonzero limit ordinal \(\lambda\le\kappa\),
\[
  P_\lambda Y=\overline{\bigcup_{\alpha<\lambda}P_\alpha Y}.
\]
\end{lemma}

\begin{proof}
Take closed subspaces \(U\subset Y\) and \(V\subset Y^*\) such that
the unit ball of \(V\) norms \(U\) and the unit ball of \(U\) norms
\(V\). The evaluation map \(J:U\to V^*\) is isometric. Its range
is closed and its annihilator in \(V^{**}=V\) is zero, so it is onto.
The operator \(P_{U,V}:Y\to U\) given by
\(P_{U,V}x=J^{-1}(f\mapsto f(x))\) is a contractive projection with
kernel \(V_\perp\).

Any two seed subspaces of density at most an infinite cardinal
\(\mu\) can be enlarged to a pair \(U,V\) with these properties
and with densities still at most \(\mu\). Adjoin norming functionals
for a dense subset of the current primal subspace, then norming
vectors for a dense subset of the enlarged dual subspace. Close
the spans after each enlargement. Norming vectors for dual
functionals exist by reflexivity. Taking the closures of the
increasing unions after countably many steps gives the desired
pair, since the two norming equalities pass to the closures.

Now enumerate a dense subset of \(Y\) as
\((y_\alpha)_{\alpha<\kappa}\). Start with the zero pair and at
successor step \(\alpha+1\) apply the enlargement to the preceding
pair and the vector \(y_\alpha\). At limits take the closures of
the two unions, which again give a mutually norming pair. The
density at each stage \(\alpha<\kappa\) is at most
\(\max\{\aleph_0,|\alpha|\}\), since a union of at most
\(|\alpha|\) preceding dense sets has that bound. Define
\(P_\alpha=P_{U_\alpha,V_\alpha}\). The ranges increase and the
kernels decrease, which gives the two projection identities.
The limit formula follows from the construction and
\(U_\kappa=Y\) by density of the enumerated set.
\end{proof}

For the projections of Lemma~\ref{lem:pri}, set
\[
  Y_\alpha=(P_{\alpha+1}-P_\alpha)Y,\qquad \alpha<\kappa .
\]
Each \(Y_\alpha\) is the range of a bounded projection and is a
closed subspace of \(P_{\alpha+1}Y\). If \(\alpha_1<\dots<\alpha_m\)
and \(u_j\in Y_{\alpha_j}\), then \(P_{\alpha_j}\) fixes the previous
partial sum and annihilates \(u_j\). Repeated use of
\eqref{eq:one-sided} therefore yields
\begin{equation}\label{eq:block-upper}
  \left\|\sum_{j=1}^m u_j\right\|^p
     \le D\sum_{j=1}^m\|u_j\|^p .
\end{equation}
The difference spaces span a dense subspace of \(Y\). Indeed,
\(P_{\alpha+1}Y=P_\alpha Y+Y_\alpha\), and at a limit ordinal
one takes the norm closure of the preceding ranges.

\section{Proof of the main theorem}\label{sec:proof}

\begin{proof}[Proof of Theorem~\ref{thm:generation}]
Fix \(p\). We use induction on density to show that every space
with a modulus of smoothness of power type \(p\) is \(\ell_p\)-generated.
The zero space is immediate. For nonzero separable \(Y\), take a sequence
\((y_n)\) in its unit ball with dense linear span. The map
\[
  (a_n)\longmapsto \sum_{n=1}^{\infty}2^{-n}a_ny_n
\]
is bounded by H\"older's inequality and has dense range.

Let \(Y\) be nonseparable with \(\dens Y=|\kappa|\), and assume
the assertion at all smaller densities. The norm of \(Y\) is
uniformly smooth, so the dual norm is uniformly convex and \(Y\)
is reflexive. Use the projections from Lemma~\ref{lem:pri}. For
each \(\alpha<\kappa\), the space \(Y_\alpha\) inherits the
power-type estimate, and
\(\dens Y_\alpha\le \max\{\aleph_0,|\alpha+1|\}<|\kappa|\).
By the induction hypothesis, choose a dense-range map
\(S_\alpha:\ell_p(\Delta_\alpha)\to Y_\alpha\) with
\(\|S_\alpha\|\le 1\). Such a normalization is possible by a
positive scalar multiplication. We use the zero map from the
zero space when \(Y_\alpha=\{0\}\).

On the finite block sums define
\[
  S(a)=\sum_{\alpha<\kappa} S_\alpha a_\alpha .
\]
The disjoint union \(\Delta=\bigsqcup_{\alpha<\kappa}\Delta_\alpha\)
identifies the \(\ell_p\)-sum of the domains with \(\ell_p(\Delta)\).
By \eqref{eq:block-upper},
\[
  \|S(a)\|^p\le D\sum_{\alpha<\kappa}\|S_\alpha a_\alpha\|^p
               \le D\|a\|_{\ell_p(\Delta)}^p .
\]
Thus \(S\) extends to \(\ell_p(\Delta)\). Its range closure contains
every \(Y_\alpha\) and hence equals \(Y\). The density bound on the
individual ranges and the closure at limits apply equally to
successor and limit cardinals. This completes the induction.

In (i), the map \(T:Y\to X\) remains bounded with dense range
after the equivalent renorming of \(Y\). Composing \(T\) with
the map just constructed gives a bounded dense-range operator
\(R:\ell_p(\Delta)\to X\). To obtain the claimed cardinality, take
a dense set in \(X\) of size \(\dens X\). For each of its points,
choose a sequence of finitely supported vectors whose images
under \(R\) converge to that point. The union of their supports
has size at most \(\dens X\), and the restriction of \(R\) to
this union still has dense range. If necessary, add coordinates
on which \(R\) is zero so that the index set has exactly the
stated cardinality.

The converse follows from the smoothness of the canonical norm
of \(\ell_p(\Gamma)\). For real
or complex scalars \(a,b\), concavity of \(s^{p/2}\) and the
scalar parallelogram identity give
\[
  \tfrac12(|a+b|^p+|a-b|^p)
  \le (|a|^2+|b|^2)^{p/2}\le |a|^p+|b|^p .
\]
Summing over the coordinates and taking \(p\)-th roots yields
\(\rho_{\ell_p(\Gamma)}(t)\le(1+t^p)^{1/p}-1\le t^p/p\).
\end{proof}

\section{Consequences}\label{sec:consequences}

The case of a superreflexive generator and the case of an
\(\ell_2\)-sum of such generators extend
\cite[Theorem~3(ii)--(iii)]{FGHZ} to all density characters.

\begin{corollary}\label{cor:superreflexive}
A Banach space is generated by a superreflexive space if and
only if it is generated by \(\ell_p(\Gamma)\) for some
\(1<p\le 2\) and some index set \(\Gamma\).
\end{corollary}

\begin{proof}
By Pisier's renorming theorem \cite[Chapter~11]{Pisier}, a
superreflexive space admits an equivalent norm whose modulus of
smoothness has power type \(p\) for some \(1<p\le 2\).
Theorem~\ref{thm:generation} gives the implication. The converse
follows from the superreflexivity of \(\ell_p(\Gamma)\).
\end{proof}

\begin{corollary}\label{cor:sum}
A Banach space is generated by an \(\ell_2\)-sum of superreflexive
spaces if and only if it is generated by a space of the form
\[
  \left(\bigoplus_{n=2}^{\infty}\ell_{1+1/n}(\Gamma_n)\right)_{\ell_2}
\]
for suitable index sets \(\Gamma_n\).
\end{corollary}

\begin{proof}
Write the given source space as
\((\bigoplus_{\lambda\in\Lambda}Y_\lambda)_{\ell_2}\).
For each \(\lambda\), use Corollary~\ref{cor:superreflexive}
to choose a generating exponent \(p_\lambda>1\). Assign it to
one group \(n\ge 2\) with \(q_n=1+1/n\le p_\lambda\).
The inclusion of \(\ell_{q_n}\) into \(\ell_{p_\lambda}\) on
the same index set has dense range. We thus obtain for each
\(\lambda\) in group \(n\) a contractive dense-range map
\(V_\lambda:\ell_{q_n}(\Gamma_\lambda')\to Y_\lambda\).
Set \(\Gamma_n=\bigsqcup_{\lambda\text{ in group }n}\Gamma_\lambda'\).
The maps \(V_\lambda\) together define a contraction from
\(\ell_{q_n}(\Gamma_n)\) to the corresponding \(\ell_2\)-sum
because \(q_n\le 2\). Finite block sums show that its range is
dense. Taking the \(\ell_2\)-sum of these maps and composing
with the original generating operator gives the desired
generation. Conversely, each \(\ell_{1+1/n}(\Gamma_n)\) is
superreflexive.
\end{proof}

For a Lebesgue space the exponent in Theorem~\ref{thm:generation}
is determined by its usual norm. The quadratic case below includes
the finite-measure question in \cite[Problem~1.3]{Rodriguez}.

\begin{corollary}\label{cor:lebesgue}
Let \(1<r<\infty\) and let \(\mu\) be a measure. Every closed
subspace of \(L^r(\mu)\) is generated by
\(\ell_{\min\{r,2\}}(\Gamma)\) for some \(\Gamma\). In
particular, for \(r\ge 2\), every such subspace is Hilbert
generated.
\end{corollary}

\begin{proof}
For \(r\le 2\), integration of the scalar inequality used in the
proof of Theorem~\ref{thm:generation} gives smoothness of power
type \(r\). For \(r>2\), the second derivative of \(z\mapsto|z|^r\)
has norm at most \(r(r-1)|z|^{r-2}\), with complex scalars
viewed as vectors in \(\mathbb R^2\). Hence
\[
  |a+b|^r+|a-b|^r-2|a|^r
  \le r(r-1)(|a|+|b|)^{r-2}|b|^2 .
\]
For \(f,g\in S_{L^r(\mu)}\) and \(0<t\le 1\), apply this estimate
with \(a=f\) and \(b=tg\) and integrate. H\"older's inequality gives
\[
  \int (|f|+t|g|)^{r-2}|g|^2\,d\mu
  \le \big\||f|+t|g|\big\|_r^{r-2}\|g\|_r^2
  \le 2^{r-2}.
\]
Taking the \(r\)-th root of the resulting bound on the mean of
\(\|f+tg\|_r^r\) and \(\|f-tg\|_r^r\) gives
\(\rho_{L^r(\mu)}(t)=O(t^2)\). Closed subspaces inherit these
estimates, and Theorem~\ref{thm:generation} applies.
\end{proof}

The exponent in the main theorem cannot be increased uniformly.
The obstruction already arises for the sequence spaces themselves,
as in \cite[Lemma~6]{FGHZ}. The following matrix argument gives
the corresponding separable-range statement.

\begin{proposition}\label{prop:sharp}
If \(1<p<r<\infty\), then every bounded linear operator
\(A:\ell_r(\Lambda)\to\ell_p(\Gamma)\) has separable range.
In particular, for uncountable \(\Gamma\) the space
\(\ell_p(\Gamma)\) is not generated by any \(\ell_r(\Lambda)\).
\end{proposition}

\begin{proof}
Write \(a_{\gamma\lambda}=(Ae_\lambda)(\gamma)\). Each column
has countable support, and so does each row, since the row is a
bounded linear functional on \(\ell_r(\Lambda)\). If
uncountably many rows were nonzero, there would be \(\varepsilon>0\)
and an uncountable set of rows \(H\) such that for each
\(\gamma\in H\) one can choose \(\lambda(\gamma)\) with
\(|a_{\gamma,\lambda(\gamma)}|>\varepsilon\).

Choose \(\gamma_j\in H\) recursively. At each step discard the
supports of previously selected columns and discard the rows
whose chosen columns belong to the support of a previously
selected row. This removes only countably many rows, by the two
countability properties above. With \(\lambda_j=\lambda(\gamma_j)\)
the resulting matrix on the selected rows and columns is
diagonal, with distinct column indices. Consequently
\[
   \varepsilon N^{1/p}\le
   \left\|A\sum_{j=1}^N e_{\lambda_j}\right\|_p
   \le \|A\|N^{1/r}\qquad(N\ge 1),
\]
a contradiction. The range of \(A\) therefore lies in the
sequence space on the countably many nonzero rows.
\end{proof}

\section*{Acknowledgements}

The authors thank the Functional Analysis Research Team at the School of Mathematics and Statistics of Anqing Normal University for discussing the manuscript
and helping to resolve difficulties in the proof.

\section*{Declarations}

\noindent\textbf{Funding}\enspace No funding was received for this work.

\smallskip
\noindent\textbf{Competing interests}\enspace
The authors declare that they have no competing interests.

\smallskip
\noindent\textbf{Data availability}\enspace
No datasets were generated or analysed in this study.


\begin{thebibliography}{9}

\bibitem{AmirLindenstrauss}
D. Amir and J. Lindenstrauss,
The structure of weakly compact sets in Banach spaces,
\emph{Ann. of Math. (2)} \textbf{88} (1968), 35--46.

\bibitem{FGZ-flat}
M. Fabian, A. Gonz\'alez and V. Zizler,
Flat sets, \(\ell_p\)-generating and fixing \(c_0\) in the nonseparable setting,
\emph{J. Aust. Math. Soc.} \textbf{87} (2009), 197--210.

\bibitem{FGHZ}
M. Fabian, G. Godefroy, P. H\'ajek and V. Zizler,
Hilbert-generated spaces,
\emph{J. Funct. Anal.} \textbf{200} (2003), 301--323.

\bibitem{FGMZ}
M. Fabian, G. Godefroy, V. Montesinos and V. Zizler,
Inner characterizations of weakly compactly generated Banach spaces
and their relatives,
\emph{J. Math. Anal. Appl.} \textbf{297} (2004), 419--455.

\bibitem{FMZ-survey}
M. Fabian, V. Montesinos and V. Zizler,
Smoothness in Banach spaces. Selected problems,
\emph{Rev. R. Acad. Cienc. Exactas F\'is. Nat. Ser. A Mat. RACSAM}
\textbf{100} (2006), 101--125.

\bibitem{GMZ}
A. J. Guirao, V. Montesinos and V. Zizler,
\emph{Open Problems in the Geometry and Analysis of Banach Spaces},
Springer, Cham, 2016.

\bibitem{Nordlander}
G. Nordlander,
The modulus of convexity in normed linear spaces,
\emph{Ark. Mat.} \textbf{4} (1960), 15--17.

\bibitem{Pisier}
G. Pisier,
\emph{Martingales in Banach Spaces},
Cambridge Studies in Advanced Mathematics, vol.~155,
Cambridge University Press, Cambridge, 2016.

\bibitem{Rodriguez}
J. Rodr\'iguez,
Open problems in Banach spaces and measure theory,
\emph{Quaest. Math.} \textbf{40} (2017), 811--832.

\end{thebibliography}
\end{document}